\documentclass[12pt]{amsart}
\usepackage{amsmath,amsfonts,amssymb} \usepackage{color}
\usepackage[all,cmtip]{xy}
\usepackage{graphicx}
 \usepackage{youngtab}

\newcommand{\corr}[1]{\langle {#1} \rangle}

  \newcommand{\cF}{\mathcal{F}} 
 
\newcommand{\bt}{{\bf t}}

 \newcommand{\bZ}{\mathbb{Z}}
 \newcommand{\bC}{\mathbb{C}}
 \newcommand{\pd}{\partial}
\newcommand{\Mbar}{\overline{\mathcal M}}
\newcommand{\bx}{\mathbf{x}}
\newcommand{\vac}{|0\rangle}

  \DeclareMathOperator{\res}{res}  \DeclareMathOperator{\Tr}{Tr}
  
\DeclareMathOperator{\Span}{span}

\newcommand{\be}{\begin{equation}}
\newcommand{\ee}{\end{equation}}
\newcommand{\bea}{\begin{eqnarray}}
\newcommand{\eea}{\end{eqnarray}}
\newcommand{\ben}{\begin{eqnarray*}}
\newcommand{\een}{\end{eqnarray*}}

\newcommand{\half}{\frac{1}{2}}

\newtheorem{cor}{Corollary}[section]

 \newtheorem{thm}[cor]{Theorem}
\theoremstyle{remark}

\definecolor{A}{rgb}{.75,1,.75}

\definecolor{yellow}{rgb}{1,1,0}
\definecolor{orange}{rgb}{1,.7,0}
\definecolor{red}{rgb}{1,0,0}
\definecolor{white}{rgb}{1,1,1}

 \makeindex
\begin{document}
\title[Witten-Kontsevich Tau-Function]
{Explicit Formula for Witten-Kontsevich Tau-Function}

\author{Jian Zhou}
\address{Department of Mathematical Sciences\\Tsinghua University\\Beijng, 100084, China}
\email{jzhou@math.tsinghua.edu.cn}

\begin{abstract}
We present an explicit formula for Witten-Kontsevich tau-function.
\end{abstract}

\maketitle

\section{Introduction}

By the famous Witten Conjecture/Kontsevich Theorem \cite{Witten1, Kon},
the generating series of intersection numbers of $\psi$-classes on moduli spaces  of algebraic curves
is a tau-function of the KdV hierarchy.
Furthermore,
as Witten \cite{Witten1} pointed out,
this fact together with the string equation completely determines
this tau-function.
Finding the explicit formula of this Witten-Kontsevich tau-function has remained an open problem
for a long time.
We will solve this problem in this paper.

In Kyoto school's approach to the integrable hierarchy \cite{MJD},
there are three pictures to describe the tau-function of the KdV hierarchy:
(a) as an element of the Sato infinite-dimensioanl Grassmannian;
(b) as a vector in the bosonic Fock space, i.e. the space of symmetric functions;
(c) as a vector in the fermionic Fock space, i.e., the semi-infinite wedge space.
Hence one can describe the Witten-Kontsevich tau-function from each of these pictures.
For the first picture,
it was given by Kac-Schwarz \cite{Kac-Schwarz} using the asymptotic expansion of the Airy function.
For the second picture,
Alexandrov \cite{Alexandrov} recently gave a recursive solution in terms of a cut-and-join type differential operator:
\be
Z_{WK} = e^{\hat{W}} 1= 1 + \sum_{n \geq 1} \frac{1}{n!} \hat{W}^n 1,
\ee
where the operator $\hat{W}$ is defined by
\ben
\hat{W}
& = &\frac{2}{3} \sum_{k,l \geq 0}(k+\frac{1}{2})u_k \cdot (l+\frac{1}{2}) u_l \cdot \frac{\pd}{\pd u_{k+l-1}} \\
& + & \frac{\lambda^2}{12} \sum_{k, l \geq 0} (k+l+\frac{5}{2}) u_{k+l+2} \frac{\pd^2}{\pd u_k\pd u_l}
+ \frac{1}{\lambda^2} \frac{u_0^3}{3!} + \frac{u_1}{16}.
\een
Inspired by these results,
we study the Witten-Kontsevich tau-function in the third picture
and obtain the following result:
The Witten-Kontsevich tau-function is a Bogoliubov transformation:
\be
Z_{WK} = e^A \vac,  \quad \quad \quad
A =  \sum_{m,n \geq 0} A_{m,n} \psi_{-m-\frac{1}{2}}\psi_{-n-\frac{1}{2}}^*,
\ee
where $A_{m,n} = 0$ if $m+n \not\equiv -1 \pmod{3}$ and
\ben
&& A_{3m-1,3n} = A_{3m-3, 3n+2}
= (-1)^n \biggl( -\frac{\sqrt{-2}}{144}\biggr)^{m+n} \frac{(6m+1)!!}{(2(m+n))!} \\
&& \qquad\qquad \cdot \prod_{j=0}^{n-1} (m+j) \cdot \prod_{j=1}^{n} (2m+2j-1) \cdot (B_{n} (m) +\frac{b_n}{6m+1}), \\
&& A_{3m-2,3n+1}
= (-1)^{n+1} \biggl( -\frac{\sqrt{-2}}{144}\biggr)^{m+n} \frac{(6m+1)!!}{(2(m+n))!} \\
&& \qquad\qquad  \cdot \prod_{j=0}^{n-1} (m+j) \cdot \prod_{j=1}^{n} (2m+2j-1) \cdot (B_{n}(m) +\frac{b_n}{6m-1}),
\een
where  $B_n(m)$ is a polynomial in $m$ of degree $n-1$ defined by:
\be
B_n(x) = \frac{1}{6} \sum_{j=1}^{n} 108^{j} b_{n-j} \cdot (x+n)_{[j-1]},
\ee
where
\be
(a)_{[j]} = \begin{cases}
1, & j = 0, \\
a(a-1) \cdots (a-j+1), & j > 0,
\end{cases}
\ee
and $b_n$ is a constant depending on $n$ defined by:
\be
b_{n} = \frac{2^n \cdot (6n+1)!!}{(2n)!}.
\ee
After using the boson-fermion correspondence,
our result can be translated into the bosonic picture
to give an explicit formula for the Witten-Kontsevich tau-function.
and it takes a very simple form.
The Witten-Konsevich tau-function is given in suitable coordinates as follows:
\be \label{eqn:Main}
Z_{WK} = \sum_{|\mu| \equiv 0 \pmod{3}} A_\mu \cdot s_\mu,
\ee
where $s_\mu$ is the Schur function indexed by a partition $\mu$,
and $A_\mu$ is a specialization of $s_\mu$ given as follows.
Let $\mu = (m_1, \dots, m_k|n_1, \dots, n_k)$ in Frobenius notation,
then
\be
A_\mu =(-1)^{n_1+ \cdots +n_k} \cdot  \det (A_{m_i,n_j})_{1 \leq i, j \leq k}.
\ee
Originally we have attempted to derive the fermionic operator
from Alexandrov's result by boson-fermion correspondence.
This is a natural approach but unfortunate we did not succeed.

Our result solves a special case of a remarkable conjecture of Aganagic-Dijkgraaf-Klemm-Mari\~no-Vafa \cite{ADKMV}.
In joint work with Fusheng Deng \cite{Deng-Zhou1, Deng-Zhou2},
we have established some other cases of this conjecture.
In \cite{Zhou1} we have generalized Alexandrov's result to Witten's r-spin intersection numbers
(this corresponds to simple singularities of type A)
and further generalization to the simple singularities of types D and E can also be made \cite{Zhou2}.
In a work in progress \cite{Zhou3},
we will generalize our result to these cases.
Working on the general $r$-spin case has provided with us some important insight on the original
Witten-Kontsevich tau-function,
which corresponds to the $r=2$ case.
In particular,
a choice of coordinates that make the idea of Bogoliubov transformation  work
comes from the consideration of the $r$-spin case in \cite{Zhou1}.

The rest of this paper is arranged as follows.
In Section 2 we review some preliminary backgrounds on partitions, symmetric functions,  
and boson-fermion correspondence.
In Section 3 we recall Witten Conjecture/Kontsevich Theorem, DVV Virasoro constraints,
and Alxexandrov's bosonic operator formula for the Witten-Kontsevich tau-function.
In Section 4 we present our main result,
whose proof is presented in the last two Sections.

\section{Preliminaries on Boson-Fermion Correspondence}

\subsection{Partitions}
A partition $\mu$ of a nonnegative integral number $n$ is a decreasing finite sequence
of nonnegative integers $\mu_1\geq\cdots \geq\mu_l>0$,
such that $|\mu|:= \mu_1 + \cdots + \mu_l = n$.
It is very useful to graphically represent a partition by its Young diagram,
e.g., $\mu=(3,2)$ is represented by:
$$\yng(3,2)$$
This leads to many natural definitions.
First of all,
by transposing the Young diagram one can define the conjugate $\mu^t$ of $\mu$,
e.g., for $\mu =(3,2)$, $|mu^t=(2,2,1)$:
$$\yng(3,2) \leftrightarrows \yng(2,2,1)$$
Secondly,
assume the Young diagram of $\mu$ has $k$ boxes in the diagonal.
Define $m_i = \mu_i - i$ and $n_i = \mu^t_i - i$ for $i = 1, \cdots , k$,
then it is clear that $m_1> \cdots > m_k \geq 0$ and  $n_1> \cdots > n_k \geq 0$.
The partition $\mu$ is completely determined by the numbers $m_i , n_i$.
One can denote the partition $\mu$ by $(m_1, \dots , m_k | n_1, \dots , n_k)$,
this is called the Frobenius notation.
A partition of the form $(m|n)$ in Frobenius notation is called a hook partition,
for example:
$$\yng(4,1,1)$$

\subsection{Schur functions}
Let $\Lambda$ be the space of symmetric functions in $\bx = (x_1, x_2, \dots)$.
For a partition $\mu$, let $s_\mu:=s_\mu(\bx)$ be the Schur function in $\Lambda$.
If we write $\mu = (m_1, \cdots , m_k | n_1, \cdots , n_k)$ in Frobenius notation, then there is a
determinantal formula that expresses $s_\mu$ in terms of $s_{(m|n)}$ (\cite[p. 47, Example 9]{Macdonald}):
$$s_\mu = \det(s_{(m_i|n_j)})_{1\leq i , j \leq k }.$$
where $\mu = (m_1, \cdots , m_k | n_1, \cdots , n_k)$ be a partition in Frobenius notation.

\subsection{Fermionic Fock space }

We say a sequence of half-integers $a_1>a_2> \cdots$ is admissible if both the set
$\mathbb{Z}_- + \frac{1}{2}\backslash \{a_1, a_2, \dots \}\subset \mathbb{Z}+\frac{1}{2}$
and the set  $\mathbb{Z}_- + \frac{1}{2}\backslash \{a_1, a_2, \dots\}$
are finite,
where  $\mathbb{Z}_-$ is the set of negative integers.

For an admissible sequence $a_1 > a_2 > \cdots$,
let $A = \{a_1, a_2, \dots\}$,
associate an element $\underline{A}$ in the half-infinite wedge space  as follows:
$$\underline{A} = \underline{a_1}\wedge \underline{a_2} \wedge \cdots.$$
The free fermionic Fock space $\mathcal{F}$ is defined as
$$\cF = \Span \{\underline{A}: \; A\subset \mathbb{Z}+\frac{1}{2}\; \text{is admissible} \}.$$
One can define an inner product on $\mathcal{F}$ by taking
$\{\underline{A}:\; A\subset \mathbb{Z}+\frac{1}{2}\; \text{is admissible} \}$ as an orthonormal basis.

\subsection{Charge decomposition}

There is a natural decomposition
$$\mathcal{F} = \bigoplus_{n\in \mathbb{Z}} \cF^{(n)},$$
where $\cF^{(n)} \subset \mathcal{F}$ the subspace spanned by $\underline{A}$ such that
$$|A\backslash \mathbb{Z}_- + \frac{1}{2}| - |\mathbb{Z}_- + \frac{1}{2}\backslash A| = n.$$
An operator on $\mathcal{F}$ is called charge 0 if it preserves the above decomposition.

The charge 0 subspace $\cF^{(0)}$ has a basis indexed by partitions:
\be
|\mu\rangle:= \underline{\mu_1 - \frac{1}{2}} \wedge \underline{\mu_2 - \frac{3}{2}}\wedge \cdots \wedge
 \underline{\mu_l -l + \frac{1}{2}}\wedge \underline{-l-\frac{1}{2}}\wedge \cdots,
\ee
where $\mu = (\mu_1, \cdots , \mu_l)$.
If $\mu = (m_1, \cdots , m_k | n_1, \cdots , n_k)$ in Frobenius notation, then
\be
\begin{split}
|\mu\rangle = & \underline{m_1+\frac{1}{2}}\wedge \cdots \wedge \underline{m_k+\frac{1}{2}}\wedge \underline{-\frac{1}{2}}
\wedge \underline{-\frac{3}{2}} \wedge \cdots \\
& \wedge \widehat{\underline{-n_k-\frac{1}{2}}} \wedge \cdots \wedge
\widehat{\underline{-n_1-\frac{1}{2}}} \wedge \cdots .
\end{split}
\ee
In particular,
when $\mu$ is the empty partition,
we get:
$$|0\rangle := \underline{-\frac{1}{2}}\wedge \underline{-\frac{3}{2}}\wedge \cdots \in \mathcal{F}.$$
It is called the fermionic vacuum vector.

\subsection{Creators and annihilators on $\mathcal{F}$}
For $r \in \mathbb{Z}+\frac{1}{2}$,
define operators $\psi_r$ and $\psi^*_r$ by
\begin{eqnarray*}
&\psi_{-r} (\underline{A}) =
\begin{cases}
(-1)^{k}\underline{a_1}\wedge\cdots\wedge\underline{a_k}\wedge \underline{r}\wedge\underline{a_{k+1}}\wedge\cdots, & \text{if $a_k > r > a_{k+1}$ for some $k$}, \\
0, &  \text{otherwise};
\end{cases}\\
&\psi^*_r(\underline{A}) =
\begin{cases}
(-1)^{k+1}\underline{a_1}\wedge\cdots\wedge \widehat{\underline{a_k}}\wedge\cdots, & \text{if $a_k = r$ for some $k$}, \\
0, &  \text{otherwise}.
\end{cases}
\end{eqnarray*}
The anti-commutation relations for these operators are
\begin{equation} \label{eqn:CR}
[\psi_r,\psi^*_s]_+:= \psi_r\psi^*_s + \psi^*_s\psi_r = \delta_{-r,s}id
\end{equation}
and other anti-commutation relations are zero.
It is clear that for $r > 0$,
\begin{align}
\psi_{r} \vac & = 0, & \psi_r^* \vac & = 0,
\end{align}
so the operators $\{\psi_{r}, \psi_r^*\}_{r > 0}$ are called the fermionic annihilators.
For a partition $\mu = (m_1, m_2, . . ., m_k | n_1, n_2, . . ., n_k)$, 
it is clear from the definition of $|\mu\rangle$ that
\be \label{eqn:Mu}
|\mu\rangle = (-1)^{n_1 + n_2 + \cdots + n_k}\prod_{i=1}^k \psi_{m_i+\frac{1}{2}} \psi_{-n_i-\frac{1}{2}}^*|0\rangle .
\ee
So the operators $\{\psi_{-r}, \psi_{-r}^*\}_{r > 0}$ are called the fermionic creators.

\subsection{The boson-fermion correspondence}
For any integer $n$, define an operator $\alpha_n$ on the fermionic Fock space $\mathcal{F}$ as follows:
\begin{equation*}
\alpha_n = \sum_{r\in \mathbb{Z} + \frac{1}{2}}:\psi_{-r}\psi^*_{r+n}:
\end{equation*}
Let
$\mathcal{B} = \Lambda[z , z^{-1}]$
be the bosonic Fock space, where $z$ is a formal variable.
Then the  boson-fermion correspondence is a linear isomorphism
$\Phi: \mathcal{F} \rightarrow \mathcal{B}$ given by
\begin{equation}
u\mapsto z^m \langle\underline{0}_m | e^{\sum_{n=1}^\infty \frac{p_n}{n}\alpha_n}u\rangle ,\ \ u\in \cF^{(m)}
\end{equation}
where $|\underline{0}_m\rangle = \underline{-\frac{1}{2}+m}\wedge \underline{-\frac{3}{2}+m}\wedge\cdots$.
It is clear that $\Phi$ induces an isomorphism between $\cF^{(0)}$ and $\Lambda$. Explicitly, this isomorphism is given by
\begin{equation}\label{eqn:boson-fermion}
|\mu\rangle \longleftrightarrow s_\mu.
\end{equation}

\section{Bosonic Representation of Witten-Kontsevich Tau-Function}

\subsection{Free energy of 2D topological gravity}

The correlators of the 2D topological gravity are defined as the following intersection numbers on the moduli space
$\Mbar_{g,n}$ of stable algebraic curves:
\be
\corr{\tau_{a_1} \cdots \tau_{a_n}}_g :=
\int_{\Mbar_{g,n}} \psi_1^{a_1} \cdots \psi_n^{a_n}.
\ee
Because the dimension of $\Mbar_{g,n}$ is given by the following formula:
\be
\dim \Mbar_{g,n} = 3g-3+n,
\ee
the following selection rule is satisfied by the correlators:
The correlator $\corr{\tau_{a_1} \cdots \tau_{a_n}}_g \neq 0$ only if
\be \label{eqn:SelRule}
a_1 + \cdots + a_n = 3g-3 + n.
\ee
This suggests the following grading:
\be
\deg t_a = 2 a + 1.
\ee
In the literature,
the following gradings have also been used: $\deg t_a = a + \frac{1}{2}$ or $\deg t_a = \frac{1}{3} (2a+1)$.
The reason for our choice of the grading will be made clear later in this paper.

The free energy is defined by
\be
F(\bt) = F_0(\bt) + F_1(\bt) + \cdots,
\ee
where the genus $g$ part of the free energy is defined by:
\be
F_g(\bt) = \corr{\exp \sum_{a \geq 0} t_a \tau_a }_g
= \sum_{n \geq 0} \frac{1}{n!} \sum_{a_1, \dots, a_n \geq 0} t_{a_1} \cdots t_{a_n}
\corr{ \tau_{a_1} \cdots \tau_{a_n}}_g.
\ee
It is a common practice to introduce a genus tracking parameter $\lambda$ and set
\be
F(\bt; \lambda) = \sum_{g \geq 0} \lambda^{2g-2} F_g(\bt).
\ee
For our purpose in this paper,
it is crucial to suppress the explicit dependence on $\lambda$.
However,
this will not cause any loss of information.
Indeed,
by the selection rule \eqref{eqn:SelRule},
one can recover $F(\bt; \lambda)$ from $F(\bt)$ as follows:
\be
F(\bt;\lambda) = F(\lambda^{-2/3} t_0, t_1, \dots, \lambda^{(2n-2)/3} t_n, \dots).
\ee

For the reader's convenience,
we recall the first few terms of $F_g$ here:
\ben
F_0 & = & \frac{1}{6} t_0^3 + \frac{1}{6} t_0^3t_1  + \frac{1}{6} t_0^3t_1^2 + \frac{1}{24}t_0^4t_2 + \cdots, \\
F_1 & = & \frac{1}{24} t_1 + \frac{1}{24} t_0t_2+ \frac{1}{48} t_1^2 + \frac{1}{72} t_1^3 + \frac{1}{12} t_0t_1t_2
  + \frac{1}{48} t_0^2 t_3  + \cdots, \\
F_2 & = & \frac{1}{1152} t_4  + \frac{29}{5760} t_2 t_3
  + \frac{1}{384} t_1t_4 + \frac{1}{1152} t_0t_5 + \cdots.
\een

\subsection{Partition function of 2D topological gravity}

It is defined by
\be
Z_{WK}(\bt) = \exp F(\bt).
\ee
This is also called the Witten-Kontsevich tau-function.
The following are the first few terms of $Z(\bt)$:
\ben
Z_{WK}(\bt) & = & 1+ (\frac{1}{6} t_0^3+ \frac{1}{24} t_1)+ (\frac{25}{144}t_0^3t_1+\frac{1}{24} t_0t_2
+\frac{25}{1152} t_1^2+\frac{1}{72} t_0^6) \\
& + & (\frac{7}{144} t_0^4t_2+ \frac{1225}{82944}t_1^3+\frac{1}{1152}t_4+\frac{49}{1728} t_0^6t_1 \\
& + & \frac{73}{6912} t_0^3t_1^2 + \frac{1}{576} t_1t_0t_2+\frac{1}{1296} t_0^9) + \cdots
\een

\subsection{Witten Conjecture/Kontsevich Theorem and Virasoro constraints}

Witten \cite{Witten1} conjectured that $Z_{WK}$ is a tau-function of the KdV hierarchy.
He also pointed out that together with the string equation satisfied by $Z_{WK}$,
this hierarchy of nonlinear differential equations uniquely determines $Z_{WK}$.
For a proof of this conjecture,
see Kontsevich \cite{Kon}.    
Dijkgraaf, E. Verlinde, H .Verlinde \cite{Dijkgraaf-Verlinde-Verlinde} and independently Fukuma, Kawai, Nakayama \cite{Fukuma-Kawai-Nakayama}
showed that that condition of being a solution of the KdV hierarchy satisfying the string equation
is equivalent to satisfying a sequence of linear differential equations called the Virasoro constraints.
In terms of correlators they are given by:
\be
\begin{split}
\corr{\tilde{\tau}_{a_0} \prod_{i=1}^n \tilde{\tau}_{a_i}}_g
= & \sum_{i=1}^n (2a_i+1) \corr{\tilde{\tau}_{a_0+a_i-1} \prod_{j \in [n]_i}
\tilde{\tau}_{a_j} }_g \\
+ & \half \sum_{b_1+b_2=a_0-2} \biggl( \corr{\tilde{\tau}_{b_1}\tilde{\tau}_{b_2} \prod_{i=1}^n \tilde{\tau}_{a_i}}_{g-1} \\
+ & \sum_{\substack{A_1 \coprod A_2 = [n]\\ g_1+g_2=g}}^s
\corr{\tilde{\tau}_{b_1} \prod_{i\in A_1} \tilde{\tau}_{a_i} }_{g_1}
\cdot \corr{\tilde{\tau}_{b_2} \prod_{i\in A_2} \tilde{\tau}_{a_i} }_{g_2} \biggr),
\end{split}
\ee
where $\tilde{\tau}_a =(2a+1)!! \cdot \tau_a$ and $[n]=\{1, \dots, n\}$,
$[n]_i = [n]-\{i\}$.
Here the summation $\sum_{\substack{A_1 \coprod A_2 = [n]\\ g_1+g_2=g}}^s$
is taken over pairs $(g_1,A_1)$ and $(g_2,A_2)$ that satisfy the stability conditions:
$$2g_1-1+|A_1| > 0, \; 2g_2-1+|A_2| > 0.$$
Withe the following change of coordinates:
\be
t_k = \frac{1}{2^k} (2k+1)!! u_k,
\ee
the Virasoro constraints can be written as the following sequence of differential equations:
\be \label{eqn:Virasoro}
\frac{\pd}{\pd u_{n+1}} Z_{WK} = \hat{L}_n Z_{WK}, \quad n \geq -1,
\ee
where the operators $\hat{L}_n$ are defined by:
\be
\hat{L}_n = \sum_{k=1}^\infty (k+\frac{1}{2}) u_k \frac{\pd}{\pd u_{k+n}}+ \frac{\lambda^2}{8}
\sum_{k=0}^{n-1} \frac{\pd^2}{\pd u_k\pd u_{n-k-1}} + \frac{u_0^2}{2\lambda^2} \delta_{n,-1}
+ \frac{\delta_{n,0}}{16}.
\ee

\subsection{Kontsevich model and Alexandrov's solution}
The Witten-Kontsevich tau-function can be represented by Kontsevich matrix integral \cite{Kon}:
\be
Z_{WK} =
\frac{\int [dA] \exp \big(- \lambda^{-1} \Tr \big(\frac{A^3}{3!} + \frac{\Lambda A}{2}\big)\big)}
{\int [dA] \exp \big(- \lambda^{-1} \Tr \big(\frac{A^3}{3!} \big)\big)},
\ee
with times given by the Miwa variables:
\be
u_k = \frac{\lambda}{2k+1} \Tr \frac{1}{\Lambda^{2k+1}}.
\ee
One can write
\be \label{eqn:GradDecomp}
Z_{WK} = \sum_{n=0}^\infty Z_{WK}^{(n)},
\ee
where
\be
Z_{WK}^{(n)} = \frac{(-1)^n}{n!(3!\lambda)^n}
\frac{\int [dA] (\Tr A^3)^n \exp \big(- \lambda^{-1} \Tr \big( \frac{\Lambda A}{2}\big)\big)}
{\int [dA] \exp \big(- \lambda^{-1} \Tr \big(\frac{A^3}{3!} \big)\big)}.
\ee
Alexandrov \cite{Alexandrov} noticed that if one defines $\deg u_k = \frac{2k+1}{3}$,
then one has
\be
\deg Z_{WK}^{(n)} = n.
\ee
Therefore, by \eqref{eqn:Virasoro},
\be
\hat{W} Z_{WK} = \hat{D} Z_{WK},
\ee
where
\be
\hat{D} = \frac{2}{3} \sum_{k=0}^\infty (k+\frac{1}{2}) u_k \frac{\pd}{\pd u_k}
\ee
is the Euler vector field for the grading,
and the operator
\be
\begin{split}
\hat{W} &= \frac{2}{3} \sum_{k=1}^\infty (k+\frac{1}{2}) u_k \hat{L}_{k-1} \\
& = \frac{2}{3} \sum_{k,l \geq 0}(k+\frac{1}{2})u_k \cdot (l+\frac{1}{2}) u_l \cdot \frac{\pd}{\pd u_{k+l-1}} \\
& + \frac{\lambda^2}{12} \sum_{k, l \geq 0} (k+l+\frac{5}{2}) u_{k+l+2} \frac{\pd^2}{\pd u_k\pd u_l}
+ \frac{1}{\lambda^2} \frac{u_0^3}{3!} + \frac{u_1}{16}
\end{split}
\ee
has
\be
\deg \hat{W} = 1.
\ee
Using the grading decomposition \eqref{eqn:GradDecomp},
one has
\be
\hat{W} Z_{WK}^{(n)} = (n+1) Z_{WK}^{(n+1)}, \quad n \geq 0.
\ee
Alexandrov then gets the following solution:
\be
Z_{WK} = e^{\hat{W}} 1.
\ee
We will not use this result in this paper,
but as mentioned in the Introduction,
the search for the fermionic version of this result serves as a motivation for our work.

\section{Main Theorem}

\subsection{A natural coordinate}
If we use the following change of coordinates introduced in \cite{Zhou1}:
\be
t_n = (-1)^n \sqrt{-2} T_{2n+1} \cdot \prod_{j=0}^n (j+\frac{1}{2}),
\ee
then we get
\ben
Z_{WK} & = & 1- \sqrt{-2} (\frac{1}{24} T_1^3+ \frac{1}{32} T_3) \\
& - & (\frac{25}{384}T_1^3T_3+\frac{5}{64}T_1T_5+\frac{25}{1024} T_3^2+\frac{1}{576} T_1^6) \\
& + & \sqrt{-2} (\frac{105}{4096} T_9 +\frac{1225}{98304} T_3^3
+ \frac{245}{2048} T_3T_1T_5 \\
& + & \frac{35}{512} T_1^2T_7+\frac{1225}{24576} T_1^3T_3^2 + \frac{35}{1536} T_1^4T_5 \\
& + & \frac{49}{18432} T_1^6T_3+ \frac{1}{41472} T_1^9) + \cdots
\een

Note
\be
u_n = \frac{(-1)^n\sqrt{-2}}{2} T_{2n+1},
\ee
one can rewrite the operators $L_n= - \frac{\pd}{\pd u_{n+1}} + \hat{L}_n$ and $\hat{W}$ in the $T$-coordinates as follows:
\ben
L_n & = & (-1)^{n+1} \sqrt{-2} \frac{\pd}{\pd T_{2n+3}}
+ \frac{(-1)^n}{2} \sum_{k=0}^\infty (2k+1) T_{2k+1} \frac{\pd}{\pd T_{2k+2n+1}} \\
& + & \frac{(-1)^n\lambda^2}{4}
\sum_{k=0}^{n-1} \frac{\pd^2}{\pd T_{2k+1}\pd T_{2n-2k-1}} - \frac{T_1^2}{4\lambda^2} \delta_{n,-1}
+ \frac{\delta_{n,0}}{16},
\een

\be
\begin{split}
\hat{W}
& = - \frac{\sqrt{-2} }{12} \sum_{k,l \geq 0}(2k+1) T_{2k+1} \cdot (2l+1) T_{2l+1}
\cdot \frac{\pd}{\pd T_{2k+2l-1}} \\
& - \frac{\sqrt{-2}}{24} \sum_{k, l \geq 0} (2k+2l+5) T_{2k+2l+5} \frac{\pd^2}{\pd T_{2k+1}\pd T_{2l+1}} \\
& - \frac{\sqrt{-2}}{24} T_1^3 - \frac{\sqrt{-2}}{32} T_3.
\end{split}
\ee
In particular,
\ben
L_{-1} & = & \sqrt{-2} \frac{\pd}{\pd T_1}
- \frac{1}{2} \sum_{k=1}^\infty (2k+1) T_{2k+1} \frac{\pd}{\pd T_{2k-1}} - \frac{T_1^2}{4\lambda^2},
\een
and
\ben
L_0 & = & - \sqrt{-2} \frac{\pd}{\pd T_3}
+ \frac{1}{2} \sum_{k=0}^\infty (2k+1) T_{2k+1} \frac{\pd}{\pd T_{2k+1}} + \frac{1}{16}.
\een

We define the following operator on the space $\bC[T_1, T_3, \dots]$:
\be
\gamma_n =
\begin{cases}
(-n) \cdot T_{-n},  & n < 0, \\
\frac{\pd}{\pd T_n}, & n > 0.
\end{cases}
\ee
Then we get
\be
L_n = (-1)^{n+1} \sqrt{-2} \gamma_{2n+3} + \frac{(-1)^n}{4} \sum_{a+b=n-1} :\gamma_{2a+1}\gamma_{2b+1}:
+ \frac{1}{16} \delta_{n,0},
\ee
in particular,
\be
L_{-1} = \sqrt{-2} \gamma_1 - \frac{1}{4} \sum_{a+b=-2} :\gamma_{2a+1}\gamma_{2b+1}:
\ee
and
\be
L_0 = -\sqrt{-2} \gamma_3 + \frac{1}{4} \sum_{a+b=-1} :\gamma_{2a+1}\gamma_{2b+1}: + \frac{1}{16}.
\ee

We also have
\ben
\hat{W} & = & - \frac{\sqrt{-2}}{12} \sum_{\substack{k, l \geq 0 \\ k+l>0}} \gamma_{-(2k+1)} \gamma_{-(2l+1)} \gamma_{2k+2l-1} \\
& - & \frac{\sqrt{-2}}{24} \sum_{k,l\geq 0}  \gamma_{-(2k+2l+5)}  \gamma_{2k+1} \gamma_{2l+1} - \frac{\sqrt{-2}}{24}  \gamma_{-1}^3
- \frac{\sqrt{-2}}{96} \gamma_{-3}.
\een
It can be written in a more compact form:
\be
\hat{W} = - \frac{\sqrt{-2}}{24} \sum_{k \geq 0} \gamma_{-(2k+1)}
\biggl(\sum_{a+b=k-2} : \gamma_{2a+1}\gamma_{2b+1}: + \frac{1}{4} \delta_{k,1}\biggr).
\ee

\subsection{Reformulation of Witten-Kontsevich tau-function in terms of Schur functions}

We make a further change of variables:
\be
T_{n} = \frac{1}{n} p_n.
\ee
Now we have
\ben
Z_{WK} & = & 1- \sqrt{-2} (\frac{1}{96}p_3 + \frac{1}{24} p_1^3) \\
& - & ( \frac{25}{1152} p_1^3p_3+ \frac{1}{64} p_1p_5 + \frac{25}{9216}p_3^2+ \frac{1}{576} p_1^6) \\
& + & \sqrt{-2}  (\frac{35}{12288}p_9+ \frac{1225}{2654208}p_3^3 + \frac{49}{6144} p_3p_1p_5 \\
&& + \frac{5}{512} p_1^2p_7 + \frac{1225}{221184} p_1^3p_3^2
+ \frac{7}{1536} p_1^4p_5 \\
&& + \frac{49}{55296} p_1^6 p_3+\frac{1}{41472} p_1^9) + \cdots.
\een
We understand $p_n$ as the Newton power function:
\be
p_n = x_1^n + x_2^n + \cdots,
\ee
and so each $Z^{(n)}_{WK}$ lies in $\Lambda^{3n}$,
the degree $3n$ subspace of the space $\Lambda$ of symmetric functions in $x_1, x_2, \cdots$.

Recall that for a partition $\nu = (\nu_1, \dots, \nu_l)$ of a positive integer,
$$p_\nu = p_{\nu_1} \cdots p_{\nu_l},$$
and for the empty partition $\emptyset$,
$p_\emptyset = 1$.
It is well-known that $\{p_\nu\}_\nu$ form an additive basis of $\Lambda$.
Another additive basis is provided by the Schur functions $\{ s_\nu\}_\mu$,
and they are related to each other by the Frobenius formula \cite{Macdonald}:
\be
p_\nu = \sum_\nu \chi^\mu_\nu s_\mu.
\ee
We use this formula to rewrite $Z_{WK}$ in terms of Schur functions:
\ben
Z & = & 1 - \frac{\sqrt{-2}}{96} (5 s_{(3)} + 7 s_{(2,1)} + 5 s_{(1^3)} ) \\
& - &  \frac{1}{9216} (385 s_{(6)} + 455 s_{(5,1)} + 0 \cdot s_{(4,2)}
+ 385 s_{(4,1^2)} \\
&& - 70 s_{(3,3)} - 50 s_{(3,2,1)} - 70 s_{(2,2,2)} \\
&& + 385 s_{(3,1^3)} + 455 s_{(2,1^4)} + 385s_{(1^6)} ) \\
& + &  \frac{\sqrt{-2}}{2654208} (85085 s_{(9)} + 95095 s_{(8,1)} + 0 \cdot s_{(7,2)} \\
&& + 85085 s_{(7,1^2)} - 8085 s_{(6,3)}  - 5775 s_{(6,2,1)} \\
&& + 87010 s_{(6,1^3)} + 6825 s_{(5,4)} + 0 \cdot s_{(5,3,1)} \\
&& - 6825 s_{(5,2^2)} + 0 \cdot s_{(5,2,1^2)} + 91910 s_{(5,1^4)} \\
&& + 5775 s_{(4^2,1)} + 8085 s_{(4,3,2)} + 0 \cdot s_{(4,3,1^2)} \\
&& + 0 \cdot s_{(4,2^2,1)} + 0 \cdot s_{(4,2,1^3)} + 87010 s_{(4,1^5)} \\
&& -1050  s_{(3^3)} + 8085 s_{(3^2,2,1)} -6825 s_{(3^2,1^3)} \\
&& + 5775 s_{(3,2^3)} + 0 \cdot s_{(3,2^2,1^2)} - 5775 s_{(3,2,1^4)} \\
&& + 85085 s_{(3,1^6)} + 6825 s_{(2^4,1)} - 8085 s_{(2^3,1^3)} \\
&& + 0 \cdot s_{(2^2,1^5)} + 95095 s_{(2,1^7)} + 85085 s_{(1^9)} ) + \cdots.
\een

\subsection{Fermionic formula for $W_{WK}$}
Now we use the formula for boson-fermion correspondence to rewrite $W_{WK}$ in the fermionic picture:
\ben
Z_{WK} & = & 1 - \frac{\sqrt{-2}}{96} (5 \psi_{-\frac{5}{2}}\psi^*_{-\frac{1}{2}}
- 7 \psi_{-\frac{3}{2}} \psi^*_{-\frac{3}{2}} + 5 \psi_{-\frac{1}{2}}\psi^*_{-\frac{5}{2}} ) \\
& - &  \frac{1}{9216} (385 \psi_{-\frac{11}{2}}\psi^*_{-\frac{1}{2}}  - 455 \psi_{-\frac{9}{2}} \psi^*_{-\frac{3}{2}}
+ 385 \psi_{-\frac{7}{2}} \psi^*_{-\frac{5}{2}}  \\
& + & 70 \psi_{-\frac{5}{2}}\psi^*_{-\frac{3}{2}} \psi_{-\frac{3}{2}} \psi^*_{-\frac{1}{2}}
- 50 \psi_{-\frac{5}{2}}\psi^*_{-\frac{5}{2}} \psi_{-\frac{1}{2}}\psi^*_{-\frac{1}{2}}
+ 70 \psi_{-\frac{3}{2}} \psi^*_{-\frac{5}{2}} \psi_{-\frac{1}{2}} \psi^*_{-\frac{3}{2}} \\
&  - & 385 \psi_{-\frac{5}{2}} \psi^*_{-\frac{7}{2}}  + 455 \psi_{-\frac{3}{2}} \psi^*_{-\frac{9}{2}}
- 385 \psi_{-\frac{1}{2}} \psi^*_{-\frac{11}{2}}  ) \\
& + &  \frac{\sqrt{-2}}{2654208} (85085 \psi_{-\frac{17}{2}} \psi^*_{-\frac{1}{2}}
 - 95095 \psi_{-\frac{15}{2}} \psi^*_{-\frac{3}{2}}
+ 85085 \psi_{-\frac{13}{2}} \psi^*_{-\frac{5}{2}} \\
& + & 8085 \psi_{-\frac{-11}{2}}\psi^*_{-\frac{3}{2}}\psi_{-\frac{3}{2}}\psi^*_{-\frac{1}{2}}
- 5775 s_{(6,2,1)} \\
& - & 87010 \psi_{-\frac{11}{2}} \psi^*_{-\frac{7}{2}} + 6825 s_{(5,4)}  \\
&& - 6825 s_{(5,2^2)}  + 91910  \psi_{-\frac{9}{2}} \psi^*_{-\frac{9}{2}} \\
&& + 5775 s_{(4^2,1)} + 8085 s_{(4,3,2)} \\
&&  - 87010  \psi_{-\frac{7}{2}} \psi^*_{-\frac{11}{2}} \\
&& -1050  s_{(3^3)} + 8085 s_{(3^2,2,1)} -6825 s_{(3^2,1^3)} \\
&& + 5775 s_{(3,2^3)}   - 5775 s_{(3,2,1^4)} \\
&& + 85085 \psi_{-\frac{5}{2}} \psi^*_{-\frac{13}{2}} + 6825 s_{(2^4,1)} - 8085 s_{(2^3,1^3)} \\
&&  - 95095 \psi_{-\frac{3}{2}} \psi^*_{-\frac{15}{2}} + 85085  \psi_{-\frac{1}{2}} \psi^*_{-\frac{17}{2}} ) + \cdots
\een
One can  check that:
\ben
Z & = & \exp \biggl( - \frac{\sqrt{-2}}{96} (5 \psi_{-\frac{5}{2}}\psi^*_{-\frac{1}{2}} - 7 \psi_{-\frac{3}{2}}\psi^*_{-\frac{3}{2}}
+ 5  \psi_{-\frac{1}{2}}\psi^*_{-\frac{5}{2}} ) \\
& - &  \frac{1}{9216} (385  \psi_{-\frac{11}{2}}\psi^*_{-\frac{1}{2}} - 455  \psi_{-\frac{9}{2}}\psi^*_{-\frac{3}{2}}
+ 385  \psi_{-\frac{7}{2}}\psi^*_{-\frac{5}{2}} \\
&& - 385 \psi_{-\frac{5}{2}}\psi^*_{-\frac{7}{2}} + 455  \psi_{-\frac{3}{2}}\psi^*_{-\frac{9}{2}} - 385 \psi_{-\frac{1}{2}}\psi^*_{-\frac{11}{2}} ) \\
& + &  \frac{\sqrt{-2}}{2654208} (85085  \psi_{-\frac{17}{2}}\psi^*_{-\frac{1}{2}} - 95095  \psi_{-\frac{15}{2}}\psi^*_{-\frac{3}{2}}
+ 85085  \psi_{-\frac{13}{2}}\psi^*_{-\frac{5}{2}}   \\
&& - 87010  \psi_{-\frac{11}{2}}\psi^*_{-\frac{7}{2}}  + 91910  \psi_{-\frac{9}{2}}\psi^*_{-\frac{9}{2}}  - 87010  \psi_{-\frac{7}{2}}\psi^*_{-\frac{11}{2}} \\
&& + 85085  \psi_{-\frac{5}{2}}\psi^*_{-\frac{13}{2}} - 95095 \psi_{-\frac{3}{2}}\psi^*_{-\frac{15}{2}} + 85085  \psi_{-\frac{17}{2}}\psi^*_{-\frac{1}{2}} ) + \cdots \biggr)1.
\een

Now we state our main theorem.

\begin{thm} \label{thm:Main}
In the fermionic pciture,
the Witten-Kontsevich tau-function is given by a Bogoliubov transformation:
\be
Z_{WK} = e^A \vac,
\ee
where the operator
\be
A = \sum_{m,n \geq 0} A_{m,n} \psi_{-m-\frac{1}{2}}\psi^*_{-n-\frac{1}{2}}
\ee
with the coefficients $A_{m,n}$ as given in the Introduction.
\end{thm}

By applying \eqref{eqn:Mu} and \eqref{eqn:boson-fermion},
one easily derives the formula \eqref{eqn:Main} in the Introduction.

The proof of this Theorem will occupy the next two Sections.

\section{Derivation of the Explicit Formula for $A_{m,n}$}

\subsection{Fermionic formulation of the Virasoro operators}

To make the computations easier, we make use of fields of operators and operator product expansions \cite{Kac}.
Consider the fields of fermionic operators:
\ben
&& \psi(z) = \sum_{r\in \bZ + \frac{1}{2}} \psi_r z^{-r-\frac{1}{2}}, \\
&& \psi^*(z) = \sum_{s\in \bZ + \frac{1}{2}} \psi_r^* z^{-s-\frac{1}{2}},
\een
The commutation relations
\be
[\psi_r, \psi_s]_+ =0, \quad [\psi_r^*, \psi_s^*] = 0, \quad
[\psi_r, \psi^*_s] = \delta_{r+s,0}
\ee
are equivalent to the following OPE's:
\ben
&& \psi(z) \psi^*(w) \sim :\psi(z)\psi^*(w): + \frac{1}{z-w}, \\
&& \psi^*(z) \psi(w) \sim :\psi^*(z)\psi(w): + \frac{1}{z-w},
\een
Write
\be
\gamma(z) = :\psi(z)\psi^*(z): = \sum_{n \in \bZ} \gamma_n z^{-n-1},
\ee
where
\be
\gamma_n = \sum_{\substack{r,s\in \bZ+1/2\\r+s=n }} :\psi_r\psi_s^*:.
\ee
Consider
\be
\eta(z) = \half (\gamma(z) + \gamma(-z)) = \sum_{n\in \bZ} \gamma_{2n+1} z^{-2n-2}.
\ee
By Wick's theorem,
\ben
\gamma(z)\gamma(w)
& = & :\psi(z)\psi^*(z): \cdot :\psi(w)\psi^*(w): \\
& = & :\psi(z)\psi^*(z) \psi(w)\psi^*(w): + \frac{:\psi(z)\psi^*(w):}{z-w} \\
& + & \frac{:\psi^*(z)\psi(w):}{z-w}  + \frac{1}{(z-w)^2} \\
& = & \frac{1}{(z-w)^2} + :\pd_w \psi(w)\psi^*(w):  + :\pd_w \psi^*(w)\psi(w): + \cdots,
\een
and so
\ben
&& \eta(z)\eta(w) \\
& = &  \frac{1}{4} (\gamma(z)\gamma(w) + \gamma(z) \gamma(-w) + \gamma(-z) \gamma(w) + \gamma(-z) \gamma(-w)) \\
& = & \frac{1}{4} \biggl(\frac{1}{(z-w)^2} + :\pd_w \psi(w)\psi^*(w):  + :\pd_w \psi^*(w)\psi(w): + \cdots \\
& + & \frac{1}{(z+w)^2} + :\pd_w \psi(-w)\psi^*(-w):  + :\pd_w \psi^*(-w)\psi(-w): + \cdots \\
& + & \frac{1}{(-z-w)^2} + :\pd_w \psi(w)\psi^*(w):  + :\pd_w \psi^*(w)\psi(w): + \cdots \\
& + & \frac{1}{(-z+w)^2} + :\pd_w \psi(-w)\psi^*(-w):  + :\pd_w \psi^*(-w)\psi(-w): + \cdots \biggr) \\
& = & \frac{z^2+w^2}{(z^2-w^2)^2} + \frac{1}{2} (:\pd_w\psi(w) \psi^*(w):+ :\pd_w\psi^*(w)\psi(w): \\
& + & :\pd_w\psi(-w) \psi^*(-w):+ :\pd_w\psi^*(-w)\psi(-w):  )  + \cdots.
\een
Hence
\ben
: \eta(w) \eta(w) : & = &  \frac{1}{2} \big(:\pd_w\psi(w) \psi^*(w):+ :\pd_w\psi^*(w)\psi(w): \\
& + & :\pd_w\psi(-w) \psi^*(-w):+ :\pd_w\psi^*(-w)\psi(-w): \big) \\
& = & \sum_{n \in \bZ}
\sum_{r+s=2n} (s-r) :\psi_{r}\psi^*_{s}: w^{-2n-2},
\een
It follows that
\be
\sum_{a+b=n} :\gamma_{2a+1}\gamma_{2b+1}:
=  \sum_{r+s=2n+2} (s-r) :\psi_r\psi_s^*: .
\ee

Recall the formula for $L_n$ is
\be
L_n = (-1)^{n+1} \sqrt{-2} \gamma_{2n+3} + \frac{(-1)^n}{4} \sum_{a+b=n-1} :\gamma_{2a+1}\gamma_{2b+1}:
+ \frac{1}{16} \delta_{n,0},
\ee
because we also have
\ben
\gamma_{2n+1} = \sum_{\substack{r,s\in \bZ+1/2\\r+s= 2n+1 }} :\psi_r\psi_s^*:,
\een
we have:
\be
\begin{split}
L_n = & (-1)^{n+1} \sqrt{-2} \sum_{r+s= 2n+3} :\psi_r\psi_s^*: \\
+ & \frac{(-1)^n}{4} \sum_{r+s=2n} (s-r) :\psi_r\psi_s^*:
+ \frac{1}{16} \delta_{n,0}.
\end{split}
\ee
This was also obtained by Kac and Schwarz \cite{Kac-Schwarz} in a different normalization.

\subsection{Constraints from $L_{-1}$}
Now we have
\ben
L_{-1} & = & \sqrt{-2}
\sum_{r_1+s_1= 1}
 :\psi_{r_1}\psi_{s_1}^*: - \frac{1}{4} \sum_{r_2+s_2=-2} (s_2-r_2) :\psi_{r_2}\psi_{s_2}^*: \\
& = & \sqrt{-2} \biggl(\psi_{\frac{1}{2}} \psi_{\frac{1}{2}}^*
+ \sum_{k=0}^\infty (\psi_{-k-\frac{1}{2}} \psi^*_{k+\frac{3}{2}}
- \psi^*_{-k-\frac{1}{2}} \psi_{k+\frac{3}{2}} ) \biggr)\\
& - &  \frac{1}{4}
\biggl(\psi_{-\frac{3}{2}} \psi^*_{-\frac{1}{2}} - \psi_{-\frac{1}{2}}\psi^*_{-\frac{3}{2}}
+ \sum_{l=0}^\infty (2l+3) (\psi_{-l-\frac{5}{2}} \psi^*_{l+\frac{1}{2}} + \psi^*_{-l-\frac{5}{2}} \psi_{l+\frac{1}{2}} )
\biggr),
\een
Under the assumption that:
\be
Z = \exp(\sum_{m,n \geq 0} A_{m,n} \psi_{-m-\frac{1}{2}} \psi^*_{-n-\frac{1}{2}}) \vac,
\ee
where $A_{m,n} \neq 0$ only if $m, n \geq 0$ and $m+n \equiv -1 \pmod{3}$,
we have
\ben
e^{-A} L_{-1} Z  & = & \biggl[\sqrt{-2} \biggl(-\sum_{m, n \geq 0} A_{m,0} A_{0,n} \psi_{-m-\frac{1}{2}} \psi^*_{-n-\frac{1}{2}} \\
& + &\sum_{k,l \geq 0} ( \psi_{-k-\frac{1}{2}} A_{k+1,l} \psi^*_{-l-\frac{1}{2}}
+ \psi^*_{-k-\frac{1}{2}} A_{l,k+1} \psi_{-l-\frac{1}{2}} ) \biggr)\\
& - &  \frac{1}{4}\biggl(\psi_{-\frac{3}{2}} \psi^*_{-\frac{1}{2}} - \psi_{-\frac{1}{2}}\psi^*_{-\frac{3}{2}} \\
& + & \sum_{k,l \geq 0} (2l+3) (\psi_{-l-\frac{5}{2}} A_{l,k} \psi^*_{-k-\frac{1}{2}}
- \psi^*_{-l-\frac{5}{2}} A_{k,l} \psi_{-k-\frac{1}{2}} ) \biggr) \biggr] \vac.
\een
By comparing the coefficients of $\psi_{-m-\frac{1}{2}}\psi^*_{-n-\frac{1}{2}}$,
we get the following  equations:
\be
\begin{split}
&  \sqrt{-2} \big(- A_{m,0} A_{0,n} + A_{m+1,n} - A_{m,n+1} \big) \\
- &  \frac{1}{4}\biggl( \delta_{m,1}\delta_{n,0} - \delta_{m,0}\delta_{n,1}  \\
+ &   (2m-1) A_{m-2,n}
+ (2n-1) A_{m,n-2}   ) \biggr) = 0.
\end{split}
\ee
This can be converted into the following  recursion relations:
\be \label{eqn:L-1Recursion}
\begin{split}
A_{m,n+1} = &   \big(A_{m+1,n}  - A_{m,0} A_{0,n}  \big)
- \frac{1}{4\sqrt{-2}}\biggl( \delta_{m,1}\delta_{n,0} - \delta_{m,0}\delta_{n,1}  \\
+ &   (2m-1) A_{m-2,n} + (2n-1) A_{m,n-2}   ) \biggr).
\end{split}
\ee
We will make an extra assumption on the initial values:
\be \label{eqn:InitVal}
A_{0, 3m-1} = (-1)^{3m-1} \biggl(- \frac{\sqrt{-2}}{144} \biggr)^m \frac{(6m-1)!!}{(2m)!},.
\ee
This will uniquely determine all $A_{m,n}$'s.
One can also convert it into the following recursion relations:
\be
\begin{split}
A_{m+1,n}  & = \big( A_{m,n+1}  + A_{m,0} A_{0,n} \big)
+ \frac{1}{4\sqrt{-2}} \biggl( \delta_{m,1}\delta_{n,0} - \delta_{m,0}\delta_{n,1}  \\
& +  (2m-1) A_{m-2,n} + (2n-1) A_{m,n-2}    \biggr) .
\end{split}
\ee
This will also determine all $A_{m,n}$'s, together with the initial value:
\be
A_{3m-1,0} = \biggl(- \frac{\sqrt{-2}}{144} \biggr)^m \frac{(6m-1)!!}{(2m)!}.
\ee
One then notice that
$F_n(m) = A_{m,n}$ and $G_n(m) = (-1)^nA_{n,m}$ satisfying the same recursion relations in $n$ and the same initial values for $n= 0$,
hence we have
\be
A_{m,n} = (-1)^nA_{n,m}.
\ee

\subsection{Ansatz for $A_{m,n}$}

We separate the recursion relations \eqref{eqn:L-1Recursion} into the following three cases:
\be
\begin{split}
& A_{3m-1, 3n} =  A_{3m,3n-1}   -  A_{3m-1,0}A_{0,3n-1}  \\
& \quad\quad - \frac{6m-3}{4 \sqrt{-2} }   A_{3m-3,3n-1}
- \frac{6n-3}{4\sqrt{-2}} A_{3m-1,3n-3},
\end{split}
\ee
\be
\begin{split}
& A_{3m-2, 3n+1} = A_{3m-1,3n}     \\
& \quad\quad - \frac{6m-5}{4 \sqrt{-2} }   A_{3m-4,3n}
- \frac{6n-1}{4\sqrt{-2}} A_{3m-2,3n-2},
\end{split}
\ee
\be
\begin{split}
& A_{3m-3, 3n+2} =  A_{3m-2,3n+1} \\
& \quad\quad - \frac{6m-7}{4\sqrt{-2}}A_{3m-5,3n+1}
- \frac{6n+1}{4\sqrt{-2}} A_{3m-3, 3n-1}.
\end{split}
\ee
The following are the first few  examples that we obtain:
\ben
&& A_{3m-1,0} = \biggl(- \frac{\sqrt{-2}}{144} \biggr)^m \frac{(6m+1)!!}{(2m)!} \cdot \frac{1}{6m+1}, \\
&& A_{3m-2,1} = - \biggl(- \frac{\sqrt{-2}}{144} \biggr)^m \frac{(6m+1)!!}{(2m)!} \cdot \frac{1}{6m-1}, \\
&& A_{3m-3,2} = \biggl(- \frac{\sqrt{-2}}{144} \biggr)^m \frac{(6m+1)!!}{(2m)!} \cdot \frac{1}{6m+1}, \\
\een
\ben
&& A_{3m-1,3} = - \biggl(- \frac{\sqrt{-2}}{144} \biggr)^{m+1} \frac{(6m+1)!!}{(2(m+1))!}
m(2m+1) (18 + \frac{105}{6m+1}), \\
&& A_{3m-2,4} = \biggl( -\frac{\sqrt{-2}}{144}\biggr)^{m+1} \frac{(6m+1)!!}{(2(m+1))!}
m (2m+1) (18 + \frac{105}{6m-1}), \\
&& A_{3m-3,5} =- \biggl( -\frac{\sqrt{-2}}{144}\biggr)^{m+1}  \frac{(6m+1)!!}{(2(m+1))!}
m(2m+1) (18 + \frac{105}{6m+1}),
\een
\ben
A_{3m-1,6} & = & \biggl( -\frac{\sqrt{-2}}{144}\biggr)^{m+2} \frac{(6m+1)!!}{(2(m+2))!}
m(m+1) \cdot (2m+3)(2m+1) \\
&&  \cdot (1944m+5778+\frac{45045}{2(6m+1)}), \\
A_{3m-2,7} & = & -\biggl( -\frac{\sqrt{-2}}{144}\biggr)^{m+2} \frac{(6m+1)!!}{(2(m+2))!}
m(m+1) \cdot (2m+1)(2m+3) \\
&&  \cdot (1944m+5778+\frac{45045}{2(6m-1)}), \\
A_{3m-3,8} & = & \biggl( -\frac{\sqrt{-2}}{144}\biggr)^{m+2} \frac{(6m+1)!!}{(2(m+2))!} m(m+1)
\cdot (2m+1)(2m+3) \\
&&  \cdot (1944m+5778+\frac{45045}{2(6m+1)})),
\een

Based on these concrete results,
we make the following assumption on the form of $A_{m,n}$: For each $n \geq 0$,
there are a polynomial
$$B_n(x) = B^{(n)}_0 x^{n-1} + B^{(n)}_1 x^{n-2} + \cdots + B^{(n)}_{n-1}$$
and a constant $b_n$ such that
\ben
&& A_{3m-1,3n} = A_{3m-3, 3n+2}
= (-1)^n \biggl( -\frac{\sqrt{-2}}{144}\biggr)^{m+n} \frac{(6m+1)!!}{(2(m+n))!} \\
&& \qquad\qquad \cdot \prod_{j=0}^{n-1} (m+j) \cdot \prod_{j=1}^{n} (2m+2j-1) \cdot (B_{n} (m) +\frac{b_n}{6m+1}), \\
&& A_{3m-2,3n+1}
= (-1)^{n+1} \biggl( -\frac{\sqrt{-2}}{144}\biggr)^{m+n} \frac{(6m+1)!!}{(2(m+n))!} \\
&& \qquad\qquad  \cdot \prod_{j=0}^{n-1} (m+j) \cdot \prod_{j=1}^{n} (2m+2j-1) \cdot (B_{n}(m) +\frac{b_n}{6m-1}).
\een
The following are some examples:
\ben
&& B_0 = 0, \\
&& B_1= 18, \\
&& B_2 = 1944x+ 5778, \\
&& B_3 = 209952 x^2  1253880x+2277477, \\
&& B_4 = 22674816x^3+226118304x^2+787643676x+1114815879, \\
&& B_5 = 2448880128x^4+36665177472x^3+207169401168x^2 \\
&& +545727699972x+2633883829515/4, \\
&& B_6 = 264479053824x^5+5546713489920x^4+46133330328000x^3 \\
&& +193184363553840x^2+424746412978761x+1828597219279695/4,
\een
and
\ben
&& b_0 = 1, \\
&& b_1 = 105, \\
&& b_2= 45045/2, \\
&& b_3 = 14549535/2, \\
&& b_4 = 25097947875/8, \\
&& b_5 = 13537833083775/8, \\
&& b_6 = 17531493843488625/16,
\een

\subsection{Recursion relations for $B_n$ and $b_n$ and their solutions}

The three types of recursion relations for $A$  lead to three types of recursion relations for $B$ as follows:
\be \label{eqn:Rec1}
\begin{split}
& B_n(x) + \frac{b_{n}}{6x+1} \\
= & - 3 (B_{n-1}(x+1) + \frac{b_{n-1}}{6x+7}) \cdot \frac{(6x+5)(6x+7)}{x} \\
+ & \frac{2^n(x+n)}{6x+1}\cdot \frac{(6n-1)!!}{(2n)! \cdot  x}
+ 216 (x+n)  \cdot (B_{n-1}(x) + \frac{b_{n-1}}{6x+1}),
\end{split}
\ee

\be \label{eqn:Rec2}
\begin{split}
& B_n(x) + \frac{b_{n}}{6x-1}
= - (B_n(x) + \frac{b_{n}}{6x+1}) \\
+ & \frac{12(x+n)(x-1) (6x-5)}{(x+n-1)(6x+1)(6x-1)}  \cdot (B_{n}(x-1) + \frac{b_{n}}{6x-5}) \\
+ & 18 (6n-1) \cdot  \frac{2(x+n)}{x+n-1}   \cdot (B_{n-1}(x) + \frac{b_{n-1}}{6x-1}),
\end{split}
\ee

\be \label{eqn:Rec3}
\begin{split}
& B_n(x) + \frac{b_{n}}{6x+1}
= - (B_n(x) + \frac{b_{n}}{6x-1}) \\
+ & \frac{12(x+n)(x-1) (6x-7)}{(x+n-1)(6x+1)(6x-1)}  \cdot (B_{n}(x-1) + \frac{b_{n}}{6x-7}) \\
+ & 18 (6n+1) \cdot  \frac{2(x+n)}{x+n-1}   \cdot (B_{n-1}(x) + \frac{b_{n-1}}{6x+1}),
\end{split}
\ee

We now find their solutions.
We take $\res_{x= -1/6}$ on both sides of \eqref{eqn:Rec1} to get:
\be
b_{n} = 36(6n-1) b_{n-1}- \frac{2^n(6n-1) \cdot (6n-1)!!}{ (2n)!}.
\ee
One easily checks that
\be \label{eqn:b}
b_{n} = \frac{2^n \cdot (6n+1)!!}{(2n)!}
\ee
is a solution.
Now we can rewrite \eqref{eqn:Rec1} as follows:
\ben
B_n(x)
& = & - 108(x+2)\cdot B_{n-1}(x+1) \\
& - & \frac{105}{x} \cdot B_{n-1}(x+1)  + \frac{18(n-1)b_{n-1}}{x} \\
& + & 216 (x+n)  \cdot B_{n-1}(x) + 18 b_{n-1} ,
\een
Subtracting \eqref{eqn:Rec3} from \eqref{eqn:Rec2} and changing $x$ to $x+1$,
one gets:
\ben
B_n(x) & = & 108 (x+2) \cdot B_{n-1}(x+1) \\
& + & 105 \frac{B_{n-1}(x+1)}{x}  - \frac{18(n-1)b_{n-1}}{x} + 18 b_{n-1}.
\een
Adding up the above two equations we get:
\be
B_n(x) = 108 (x+n) \cdot B_{n-1}(x) + 18 b_{n-1}.
\ee
From this one easily gets the following solution:
\be \label{eqn:B}
B_n(x) = \frac{1}{6} \sum_{j=1}^{n} 108^{j} b_{n-j} \cdot (x+n)_{[j-1]}.
\ee
where
\be
(a)_{[j]} = \begin{cases}
1, & j = 0, \\
a(a-1) \cdots (a-j+1), & j > 0.
\end{cases}
\ee

\section{Verifications of the Virasoro Constraints}

In this section we verify that the recursion relations from all $L_n$ are indeed satisfied for
$e^A \vac$, with the operator $A$ as derived in last section,
hence completing the proof of Theorem \ref{thm:Main}.

\subsection{Checking the equation \eqref{eqn:Rec1}}
We plug the formula \eqref{eqn:B} into \eqref{eqn:Rec1}.
The light-hand side becomes:
\ben
LHS & = & \frac{1}{6} \sum_{j=1}^{n} 108^j b_{n-j} (x+n)_{[j-1]}
+ \frac{b_n}{6x+1},
\een
the right-hand side becomes
\ben
RHS & = & -\frac{1}{2} \sum_{j=1}^{n-1} 108^j b_{n-1-j} (x+n)_{[j-1]}
\cdot \biggl(36 (x+2) + \frac{35}{x} \biggr) \\
& - & 3b_{n-1} \biggl(6+\frac{5}{x} \biggr)
+\frac{2^n(6n-1)!!}{(2n)!} \cdot \biggl( \frac{n}{x} - \frac{6n-1}{6x+1}  \biggr) \\
& + & 216 (x+n) \cdot \frac{1}{6} \sum_{j=1}^{n-1} 108^jb_{n-1-j} \cdot (x+n-1)_{[j-1]}  \\
& + & 36 b_{n-1} + \frac{36(6n-1)b_{n-1}}{6x+1}.
\een
Note
\ben
&&  \frac{15b_{n-1}}{x} + \frac{2^n(6n-1)!!}{(2n)!} \cdot  \frac{n}{x}  \\
& = & - \frac{15b_{n-1}}{x} + \frac{3(6n-1)b_{n-1}}{x} \\
& = &  \frac{18(n-1)b_{n-1}}{x}.
\een
Hence we need to show that
\ben
&& \frac{1}{2} \sum_{j=1}^{n-1} 108^j b_{n-1-j} (x+n)_{[j-1]}
\cdot \biggl(36 (x+2) + \frac{35}{x} \biggr) \\
& - &   \frac{18(n-1)b_{n-1}}{x}  =  \frac{1}{6} \sum_{j=1}^{n-1} 108^{j+1}b_{n-1-j} \cdot (x+n)_{[j]}.
\een
After multiplying by $6$ on both sides it becomes:
\ben
&& (x+2) \sum_{j=1}^{n-1} 108^{j+1} b_{n-1-j} \cdot (x+n)_{[j-1]} \\
& + & \frac{105}{x} \sum_{j=1}^{n-1} 108^j b_{n-1-j} \cdot (x+n)_{[j-1]} -  \frac{108(n-1)b_{n-1}}{x}  \\
& = &  \sum_{j=1}^{n-1} 108^{j+1}b_{n-1-j} \cdot (x+n)_{[j]}.
\een
We multiply both sides of this equation by $x$,
and note that
\ben
&& (x+n)_[j]-(x+n)_{[j-1]} = (n-1-j) \cdot (x+n)_{[j-1]},
\een
so one needs to check that:
\ben
&& 105 \sum_{j=1}^{n-1} 108^j b_{n-1-j} \cdot (x+n)_{[j-1]}   \\
& = &  108(n-1)b_{n-1}+ \sum_{j=1}^{n-2} 108^{j+1}b_{n-1-j} \cdot (n-1-j) \cdot x (x+n)_{[j-1]}.
\een
This can be done by rewriting the right-hand side in the basis $\{(x+n)_{[j]}\}_{j \geq 0}$
for $\bC[x]$.
More precisely,
the right-hand side of this equation can be rewritten as:
\ben
RHS & = &  108(n-1)b_{n-1}+ \sum_{j=1}^{n-2} 108^{j+1}b_{n-1-j} \cdot (n-1-j) \cdot (x+n)_{[j-1]} \\
&& \cdot [(x+n-j+1)-(n-j+1)] \\
& = &  108(n-1)b_{n-1}+ \sum_{j=1}^{n-2} 108^{j+1}b_{n-1-j} \cdot (n-1-j) \cdot (x+n)_{[j]} \\
& - & \sum_{j=1}^{n-2} 108^{j+1}b_{n-1-j} \cdot (n-1-j)(n-j+1) \cdot (x+n)_{[j-1]} \\
& = &  \sum_{j=1}^{n-1} 108^{j}b_{n-j} \cdot (n-j) \cdot (x+n)_{[j-1]} \\
& - & \sum_{j=1}^{n-2} 108^{j+1}b_{n-1-j} \cdot (n-1-j)(n-j+1) \cdot (x+n)_{[j-1]}.
\een
By comparing the coefficients of $(x+n)_{[j-1]}$ on both sides,
one needs to check that
\ben
&& 105  \cdot b_{n-1-j}
= b_{n-j} \cdot (n-j) - 108b_{n-1-j} \cdot (n-1-j)(n-j+1).
\een
After changing $n-j$ to $m$,
one gets an equivalent equation:
\ben
&& 105  \cdot b_{m-1}
= m b_m- 108b_{m-1} \cdot (m-1)(m+1).
\een
This last equation is equivalent to
\be
b_m = \frac{3(6m+1)(6m-1)}{m} b_{m-1}.
\ee
This is readily checked by \eqref{eqn:b}.

\subsection{Checking the equation \eqref{eqn:Rec2} and \eqref{eqn:Rec3}}
We multiply both sides of \eqref{eqn:Rec2} by $(x+n-1)(36x^2-1)$ to get:
\be
\begin{split}
& 2 (x+n-1)(36x^2-1) B_n(x) + 12 b_{n}x (x-n+1)\\
= & 12(x+n)(x-1) (6x-5)B_{n}(x-1)
+ 12 b_n (x+n)(x-1) \\
+ & 36 (6n-1) \cdot (x+n)(36x^2-1)   \cdot B_{n-1}(x) \\
+ & 36 (6n-1) \cdot (x+n)(6x+1) b_{n-1}.
\end{split}
\ee
This can be checked using similar   expansions as in last subsection.
Similarly for \eqref{eqn:Rec3}.

\subsection{Checking constraints from $L_n$ ($n \geq 0$)}
Recall
\ben
L_0 & = & - \sqrt{-2} \sum_{r_1+s_1= 3} :\psi_{r_1}\psi_{s_1}^*:
 + \frac{1}{4} \sum_{r_2+s_2=0} (s_2-r_2) :\psi_{r_2}\psi_{s_2}^*: + \frac{1}{16} \\
& = & - \sqrt{-2} \biggl(\psi_{\frac{5}{2}}\psi^*_{\frac{1}{2}} + \psi_{\frac{3}{2}} \psi^*_{\frac{3}{2}}
+ \psi_{\frac{1}{2}}\psi^*_{\frac{5}{2}}
+ \sum_{k=0}^\infty (\psi_{-k-\frac{1}{2}}\psi^*_{k+\frac{7}{2}} - \psi^*_{-k-\frac{1}{2}}\psi_{k+\frac{7}{2}} ) \biggr) \\
& + & \frac{1}{4} \sum_{l=0}^\infty (2l+1) (\psi_{-l-\frac{1}{2}}\psi^*_{l+\frac{1}{2}} + \psi^*_{-l-\frac{1}{2}}\psi_{l+\frac{1}{2}} )
+ \frac{1}{16},
\een
and so
\ben
L_0 Z& = & - e^A \biggl[\sqrt{-2} \biggl(A_{2,0} + A_{1,1} + A_{0,2} \\
& - & \sum_{m,n \geq 0} (A_{m,0}A_{2,n} + A_{m,1}A_{1,n} + A_{m,2}A_{0,n})  \psi_{-m-\frac{1}{2}}  \psi^*_{-n-\frac{1}{2}} \\
& + & \sum_{k, l \geq 0} ( \psi_{-k-\frac{1}{2}} A_{k+3,l} \psi^*_{-l-\frac{1}{2}}
+ \psi^*_{-k-\frac{1}{2}} A_{l,k+3} \psi_{-l-\frac{1}{2}} ) \biggr)\\
& - &  \frac{1}{4}  \sum_{k,l\geq 0} (2l+1) (\psi_{-l-\frac{1}{2}} A_{l,k} \psi^*_{-k-\frac{1}{2}}
- \psi^*_{-l-\frac{1}{2}} A_{k,l} \psi_{-k-\frac{1}{2}} )
- \frac{1}{16} \biggr] \vac.
\een
hence the constraints from $L_0$ are
\ben
&& \sqrt{-2} (A_{2,0} + A_{1,1} + A_{0,2}  ) - \frac{1}{16}  = 0,
\een
which is readily checked,
and
\ben
&& \sqrt{-2} \biggl(- \sum_{i=0}^2 A_{m,i}A_{2-i,n}  + A_{m+3,n} - A_{m,n+3}   \biggr) \\
& - &  \frac{1}{4}  \biggl (2m+1) A_{m,n} + (2n+1) A_{m,n} \biggr)
= 0.
\een
The latter can be divided into three cases:
\ben
A_{3m-1,3n+3} & = & A_{3m+2,3n} - A_{3m-1,0}A_{2,3n}
- \frac{6(m+n)}{4\sqrt{-2}} A_{3m-1,3n}, \\
A_{3m-2,3n+4} & = & A_{3m+1,3n+1}  - A_{3m-2,1}A_{1,3n+1}
- \frac{6(m+n)}{4\sqrt{-2}} A_{3m-2,3n+1}, \\
A_{3m-3,3n+5} & = & A_{3m,3n+2} - A_{3m-3,2}A_{0,3n+2}
- \frac{6(m+n)}{4\sqrt{-2}} A_{3m-3,3n+2}.
\een
They are equivalent to the following two relations for $B_n(x)$:
\ben
&& B_{n+1}(x) + \frac{b_{n+1}}{6x+1} \\
& = & - \frac{3(6x+5)(6x+7)}{x} \biggl( B_n(x+1) + \frac{b_n}{6x+7} \biggr) \\
& + & 216 (x+n+1)\biggl( B_n(x) + \frac{b_n}{6x+1} \biggr) \\
& + & \frac{2^{n+1}(6n+5)!!}{(2n+2)!} \cdot \frac{x+n+1}{x(6x+1)},
\een
and
\ben
&& B_{n+1}(x) + \frac{b_{n+1}}{6x-1} \\
& = & - \frac{3(6x+5)(6x+7)}{x} \biggl( B_n(x+1) + \frac{b_n}{6x+5} \biggr) \\
& + & 216 (x+n+1)\biggl( B_n(x) + \frac{b_n}{6x-1} \biggr) \\
& - & \frac{2^{n+1}(6n+7)!!}{(2n+2)!(6n+5)} \cdot \frac{x+n+1}{x(6x-1)}.
\een
They can be checked by the same method.
Similarly, one can also check the constraints for $L_1$ and $L_2$.
Because the Virasoro commutation relation
\be
[L_m, L_n] = (m-n) L_{m+n},
\ee
all $L_n$-constraints are satisfied.
This completes the proof.

\vspace{.2in}
{\em Acknowledgements}.
This research is partially supported by NSFC grant 1171174.

\end{document}